\documentclass[11pt]{article}
\usepackage{amsmath}
\usepackage{amssymb}
\usepackage{amscd}

\renewcommand{\theequation}{\arabic{section}.\arabic{equation}}

\title{Connecting  exceptional orthogonal polynomials of different kind}
\author{C. Quesne\thanks{Electronic mail: christiane.quesne@ulb.be} \\
{\small \textsl{Physique Nucl\'{e}aire Th\'{e}orique et Physique Math%
\'{e}matique, Universit\'{e} Libre de Bruxelles,}} \\
{\small \textsl{Campus de la Plaine CP229, Boulevard~du Triomphe, B-1050
Brussels, Belgium}}}
\date{ }

\begin{document}

\maketitle

\begin{abstract}
The known asymptotic relations interconnecting Jacobi, Laguerre, and Hermite classical orthogonal polynomials are generalized to the corresponding exceptional orthogonal polynomials of codimension $m$. It is proved that $X_m$-Laguerre exceptional orthogonal polynomials of type I, II, or III can be obtained as limits of $X_m$-Jacobi exceptional orthogonal polynomials of the same type. Similarly, $X_m$-Hermite exceptional orthogonal polynomials of type III can be derived from $X_m$-Jacobi or $X_m$-Laguerre ones. The quadratic transformations expressing Hermite classical orthogonal polynomials in terms of Laguerre ones is also extended to even $X_{2m}$-Hermite exceptional orthogonal polynomials.  
\end{abstract}


\section{INTRODUCTION}

Exceptional orthogonal polynomials (EOPs) are complete families of orthogonal polynomials that arise as eigenfunctions of a Sturm-Liouville eigenvalue problem \cite{gomez09}. In contrast with the classical orthogonal polynomials (COPs) of Jacobi, Laguerre, and Hermite types \cite{koekoek}, there are some gaps in the sequence of their degrees, the total number of missing ``exceptional'' degrees being known as the codimension. In addition, the corresponding differential equation contains rational coefficients instead of polynomial ones. Due to these facts, EOPs circumvent the strong limitations of Bochner's classification theorem, which characterizes Sturm-Liouville COP systems.\par
%
%
During the last few years, a lot of research activity has been devoted to the study of EOPs both from a mathematical viewpoint and for their applications in mathematical physics (for a recent account see \cite{garcia19} and references quoted therein).\par
%
%
In quantum mechanics, for instance, they have been shown to lead to some exactly solvable rational extensions of well-known quantum potentials. The $X_1$ EOPs of codimension one of Ref.~\cite{gomez09} were indeed shown to be related to the Darboux transformation in the context of shape invariant potentials in supersymmetric quantum mechanics \cite{cq08}. Soon after the introduction of the first $X_2$ EOPs of codimension two \cite{cq09}, infinite families of shape invariant potentials were introduced in relation of $X_m$ EOPs of arbitrary codimension $m$ \cite{odake09}. Multi-indexed families of $X_{m_1m_2\ldots m_k}$ EOPs, connected with multi-step Darboux transformations, were then considered \cite{gomez12, odake11}. Very recently, confluent Darboux transformations have been used to generate EOP families with an arbitrary number of real parameters \cite{garcia21}.\par
%
%
{}From a mathematical viewpoint, three main questions have been the subject of research activity in relation to EOPs. The first one is the study of the interlacing and asymptotic behaviour of their zeros (see, e.g., \cite{gomez13, kuijlaars}). The second one has to do with the recurrence relations they satisfy and which may be of two different types: either recurrence relations of order $2N+1$, where $N$ is the number of Darboux steps, with coefficients that are functions of $x$ and $n$ \cite{gomez14, odake13}, or recurrence relations of order $2m+3$, where $m$ is the codimension, with coefficients only dependent on $n$ \cite{duran, miki, odake16}. Finally, the quest for a complete classification of EOPs is a fundamental problem \cite{garcia19}, which, as shown in \cite{garcia21}, remains rather open.\par
%
%
The purpose of the present paper is to look into the counterpart of a well-known property of Jacobi, Laguerre, and Hermite COPs, namely that they are interconnected through limit relations \cite{koekoek}. As a first step in this inquiry, we plan to consider here $X_m$-Jacobi, $X_m$-Laguerre, and $X_m$-Hermite EOPs, the first two existing in three different types I, II, and III, while the last ones only belong to type III.\par
%
%
In Section~2, we list the sets of EOPs that we are going to consider, we provide their explicit expressions in terms of COPs, and we explain how the latter are related to some other definitions found in the literature. In Section~3, we present and prove limit relations connecting these EOPs. In Section~4, we briefly comment on another type of relations generalizing the quadratic transformations relating Laguerre and Hermite COPs. Finally, Section~5 contains the conclusion.\par
%
%
\section{\boldmath SETS OF $X_m$ EOPS AND THEIR EXPRESSIONS IN TERMS OF COPS}

The sets of $X_m$ EOPs that we are going to consider here mainly come from Ref.~\cite{miki} (except for type II $X_m$-Laguerre EOPs, which are not considered there), where they are given in monic form. To establish relations with other definitions where standard EOPs and COPs are used, let us remind that one goes from standard Hermite, Laguerre, and Jacobi COPs to monic ones by the following transformations
\begin{eqnarray}
  H_n(x) & \to & 2^n H_n(x), \nonumber \\
  L_n^{(\alpha)}(x) & \to & \frac{(-1)^n}{n!} L_n^{(\alpha)}(x),  \label{eq:monic} \\
  P_n^{(\alpha,\beta)}(x) & \to & \frac{\Gamma(2n+\alpha+\beta+1)}{2^n n! \Gamma(n+\alpha+\beta+1)}
      P_n^{(\alpha,\beta)}(x).  \nonumber
\end{eqnarray}
For future use, we list some relations satisfied by monic COPs in Appendix A.\par
%
%
\subsection{\boldmath $X_m$-Hermite polynomials}

It is well known \cite{gomez14, marquette} that only type III $X_m$-Hermite EOPs exist and that they are restricted to even $m$ values. They may be defined as
\begin{align}
  & \hat{H}^{({\rm III},m)}_n(x) \nonumber \\
  & \quad = \begin{cases}
    1 & \text{if $n=0$}, \\
    {\rm i}^m \left[H_m({\rm i}x) H_{n-m}(x) + {\rm i} \frac{m}{2} H_{m-1}({\rm i}x)H_{n-m-1}(x)\right]
        & \text{if $n=m+1, m+2. \ldots$},
  \end{cases} \label{eq:h}
\end{align}
with $m$ missing degrees $n=1, 2, \ldots, m$. This relation comes from \cite{miki}, where
\begin{align*}
  & \hat{H}^{({\rm III},m)}_n(x) \nonumber \\
  & \quad = \begin{cases}
    1 & \text{if $n=0$}, \\
    - \frac{1}{2}{\rm i}^m \left\{H_m({\rm i}x) (\partial_x-2x)H_{n-m-1}(x) - [\partial_x H_m({\rm i}x)] 
        H_{n-m-1}(x)\right\} & \text{if $n=m+1, m+2, \ldots$},
  \end{cases}
\end{align*}
after applying (\ref{eq:A1}) and (\ref{eq:A2}). It is connected to the previously used $y^{(m)}_n(x)$ of Ref.~\cite{marquette} by the relation
\begin{equation*}
  y^{(m)}_n(x) = \begin{cases}
     \hat{H}^{({\rm III},m)}_0(x) & \text{if $n=0$}, \\
     - 2^n \hat{H}^{({\rm III},m)}_n(x) & \text{if $n=m+1,m+2,\ldots$}.
  \end{cases}
\end{equation*}
\par
%
%
\subsection{\boldmath $X_m$-Laguerre polynomials}

$X_m$-Laguerre EOPS exist in three different types, which may be defined as
\begin{equation}
  \hat{L}^{({\rm I},m)}_n(x;\alpha) = (-1)^m \left(L_m^{(\alpha+1)}(-x) L^{(\alpha)}_n(x) - n 
  L_m^{(\alpha)}(-x) L^{(\alpha+1)}_{n-1}(x)\right), \qquad n = 0, 1, 2, \ldots,  \label{eq:l1} 
\end{equation}
\begin{align}
  \hat{L}^{({\rm II},m)}_n(x;\alpha) &= \frac{1}{\alpha+n-m} \left(nx L_m^{(-\alpha)}(x) L_{n-1}
    ^{(\alpha+1)}(x) - (m-\alpha) L_m^{(-\alpha-1)}(x) L_n^{(\alpha)}(x)\right), \nonumber \\[0.2cm] 
 & \quad   n=0,1,2,\ldots,     \label{eq:l2}
\end{align}
and
\begin{align}
  &\hat{L}^{({\rm III},m)}_n(x;\alpha) \nonumber \\
  &\quad = \begin{cases}
     1 & \text{if $n=0$}, \\
     (-1)^m \left(L_m^{(-\alpha)}(-x) L_{n-m}^{(\alpha-1)}(x) - mx L_{m-1}^{(-\alpha+1)}(-x) 
        L_{n-m-1}^{(\alpha)}(x)\right) & \text{if $n=m+1,$} \\
     &\quad \text{$m+2,\ldots$}.  
  \end{cases} \label{eq:l3}
\end{align}
For the first two types, the polynomial degrees are $n+m=m, m+1, m+2, \ldots$ (with degrees 0, 1, \ldots, $m-1$ missing), while for the third type, they are $n=0, m+1, m+2, \ldots$ (with degrees 1, 2, \ldots, $m$ missing).\par
%
%
Equation (\ref{eq:l1}) is derived from the corresponding equation of Ref.~\cite{miki}, which reads
\begin{align*}
  \hat{L}^{({\rm I},m)}_n(x;\alpha) &= (-1)^{m+1}\left\{L_m^{(\alpha)}(-x)(\partial_x-1)L_n^{(\alpha)}(x)
     - \left[\partial_x L_m^{(\alpha)}(-x)\right] L_n^{(\alpha)}(x)\right\}, \\
  &\quad n=0, 1, 2, \ldots,
\end{align*}
after using (\ref{eq:A3}) and (\ref{eq:A4}). It is directly related to the definition of type I Laguerre EOPs $L^{{\rm I},\alpha}_{m,n}(x)$, given in Eq.~(3.2) of \cite{liaw16} in terms of standard Laguerre polynomials (see also \cite{gomez13}). After transforming the latter into monic ones, one indeed gets
\begin{equation*}
  L^{{\rm I},\alpha}_{m,n}(x) = \frac{(-1)^{n-m}}{m! (n-m)!} \hat{L}^{({\rm I},m)}_{n-m}(x;\alpha-1), \qquad
  n=m, m+1, m+2, \ldots.
\end{equation*}
\par
%
%
To obtain Eq.~(\ref{eq:l2}), we started from the definition of type II Laguerre EOPs $L^{{\rm II},\alpha}_{m,n}(x)$ given in \cite{liaw16} (see also \cite{gomez13}), converted standard Laguerre COPs into monic ones and multiplied the result by an appropriate factor to get a monic EOP. As a result, we got
\begin{equation*}
  L^{{\rm II},\alpha}_{m,n}(x) = (-1)^{n-1} \frac{\alpha+1+n-2m}{m! (n-m)!} \hat{L}^{({\rm II},m}_{n-m}(x;
  \alpha+1), \qquad n=m, m+1, m+2, \ldots.
\end{equation*}
\par
%
%
{}Finally, Eq.~(\ref{eq:l3}) is derived from the corresponding equation of Ref.~\cite{miki}, namely
\begin{align}
  &\hat{L}^{({\rm III},m)}_n(x;\alpha) \nonumber \\
  &\quad = \begin{cases}
     1 & \text{if $n=0$}, \\
     (-1)^{m+1} \Bigl\{L_m^{(-\alpha)}(-x) (x\partial_x+\alpha-x) L_{n-m-1}^{(\alpha)}(x) & \\
     \quad - x \left[\partial_x L_m^{(-\alpha)}(-x)\right] L_{n-m-1}^{(\alpha)}(x)\Bigr\} & \text{if 
     $n=m+1,m+2,\ldots$},  
  \end{cases} \label{eq:l3-bis}
\end{align}
by using (\ref{eq:A3}) and (\ref{eq:A5}). After some transformations, it can be related to the definition of type III Laguerre EOPs $L^{{\rm III},\alpha}_{m,n}(x)$ given in Eq.~(5.12) of \cite{liaw16}. To show this, let us start from (\ref{eq:l3-bis}) and use (\ref{eq:A5}), thereby getting
\begin{align*}
  \hat{L}^{({\rm III},m)}_n(x;\alpha) &= 
  (-1)^m \Bigl\{L_m^{(-\alpha)}(-x)  \Bigl[L_{n-m}^{(\alpha-1)}(x) + (\alpha-x) L_{n-m-1}^{(\alpha)}(x)
       \Bigr]\\
  & \quad - L_{m+1}^{(-\alpha-1)}(-x) L_{n-m-1}^{(\alpha)}(x)\Bigr]\Bigr\}, \qquad n=m+1, m+2, \ldots.
\end{align*}
On the right-hand side of this equation, the factor between square brackets becomes
\begin{align*}
  &L_{n-m}^{(\alpha-1)}(x) + (\alpha-x) L_{n-m-1}^{(\alpha)}(x) \\
  & \quad = L_{n-m}^{(\alpha)}(x) + (n-m+\alpha-x) L_{n-m-1}^{(\alpha)}(x) \\
  & \quad = -(n-m-1) \Bigl[L_{n-m-1}^{(\alpha)}(x) + (n-m+\alpha-1) L_{n-m-2}^{(\alpha)}(x)\bigr] \\
  & \quad = -(n-m-1) x L_{n-m-2}^{(\alpha+1)}(x), 
\end{align*}
after successively using (\ref{eq:A4}), (\ref{eq:A6}), and (\ref{eq:A7}). The result reads
\begin{align}
  &\hat{L}^{({\rm III},m)}_n(x;\alpha) \nonumber \\
  &\quad = \begin{cases}
     1 & \text{if $n=0$}, \\
     (-1)^{m+1} \Bigl[(n-m-1) x L_m^{(-\alpha)}(-x) L_{n-m-2}^{(\alpha+1)}(x) & \\
     \quad + L_{m+1}^{(-\alpha-1)}(-x) L_{n-m-1}^{(\alpha)}(x)\Bigr] & \text{if 
     $n=m+1,m+2,\ldots$},  
  \end{cases} \label{eq:liaw}
\end{align}
which may then be compared with Eq.~(5.12) of \cite{liaw16} after rewriting there standard Laguerre polynomials in terms of monic ones. Hence,
\begin{equation*}
  L^{{\rm III},\alpha}_{m,n}(x) = \frac{(-1)^{n-m-1}}{m! (n-m-1)!} \hat{L}^{({\rm III},m)}_{n}(x, \alpha+1),
  \qquad n=m+1, m+2, \ldots,
\end{equation*}
and of course $L^{{\rm III},\alpha}_{m,0}(x) = \hat{L}^{({\rm III},m)}_0(x;\alpha+1) = 1$.\par
%
%
It is worth observing that in \cite{liaw16}, for studying properties of EOPs, the range of the $\alpha$ parameter has been restricted to some interval ($\alpha>0$ for type I, $\alpha>m-1$ for type II, and $0<\alpha<1$ for type III). Here, as in \cite{miki}, we consider their formal definitions without restriction on $\alpha$.\par
%
%
\subsection{\boldmath $X_m$-Jacobi polynomials}

$X_m$-Jacobi EOPs exist in three different types, which may be defined as
\begin{align}
  \hat{P}^{({\rm I},m)}_n(x;\alpha, \beta) & = \frac{1}{\beta+n-m}\Bigl\{P_m^{(\alpha,-\beta)}(x) \Bigl[
      n(1+x)P_{n-1}^{(\alpha+1,\beta+1)}(x) + \beta P_n^{(\alpha,\beta)}(x)\Bigr] \nonumber \\
  & - m(1+x) P_{m-1}^{(\alpha+1,-\beta+1)}(x) P_n^{(\alpha,\beta)}(x)\Bigr\}, \qquad n=0,1,2,\ldots,
      \label{eq:j1}
\end{align}
\begin{align}
  \hat{P}^{({\rm II},m)}_{n}(x;\alpha, \beta) &=\frac{1}{m-n-\alpha} \Bigl\{P_m^{(-\alpha,\beta)}(x) \Bigl[
     n(1-x) P_{n-1}^{(\alpha+1,\beta+1)}(x) - \alpha P_n^{(\alpha,\beta)}(x)\Bigr] \nonumber \\
  & -m(1-x) P_{m-1}^{(-\alpha+1,\beta+1)}(x) P_n^{(\alpha,\beta)}(x)\Bigr\}, \qquad n=0,1,2,\ldots,
     \label{eq:j2}
\end{align}
and
\begin{align}
  &\hat{P}^{({\rm III},m)}_n(x;\alpha,\beta) \nonumber \\
  &\quad = \begin{cases}
     1 & \text{if $n=0$}, \\
     \frac{1}{\alpha+\beta+n-2m-1} \Bigl[(\alpha+\beta+n-m-1) P_m^{(-\alpha,-\beta)}(x) 
         P_{n-m}^{(\alpha-1,\beta-1)}(x) & \\
     \quad + m(1-x^2) P_{m-1}^{(-\alpha+1,-\beta+1)}(x) P_{n-m-1}^{(\alpha,\beta)}(x)\Bigr], & \text{if 
         $n=m+1$}, \\
      & \quad \text{$m+2,\ldots$}. \\  
  \end{cases}  \label{eq:j3}
\end{align}
As in the case of $X_m$-Laguerre EOPs, $X_m$-Jacobi EOPs are of degree $n+m=m, m+1, m+2, \ldots$ (with degrees 0, 1, \ldots, $m-1$ missing) in the case of type I or II, while they are of degree $n=0, m+1, m+2, \ldots$ (with degrees 1, 2, \ldots, $m$ missing) in that of type III.\par
%
%
These relations are derived from the corresponding definitions of \cite{miki}, namely\footnote{Note that in \cite{miki}, Jacobi polynomials are denoted by $J_n^{(\alpha,\beta)}(x)$ instead of $P_n^{(\alpha,\beta)}(x)$.}
\begin{align*}
  \hat{P}^{({\rm I},m)}_n(x;\alpha, \beta) &= \frac{1}{\beta+n-m}\Bigl\{P_m^{(\alpha,-\beta)}(x) [
      (1+x)\partial_x+\beta]P_n^{(\alpha,\beta)}(x) \nonumber \\
  & - (1+x)\bigl[\partial_x P_{m}^{(\alpha,-\beta)}(x)\bigr] P_n^{(\alpha,\beta)}(x)\Bigr\}, \qquad     
      n=0,1,2,\ldots,
\end{align*}
\begin{align*}
  \hat{P}^{({\rm II},m)}_{n}(x;\alpha, \beta) &= \frac{1}{m-n-\alpha} \Bigl\{P_m^{(-\alpha,\beta)}(x) [
      (1-x)\partial_x-\alpha] P_{n}^{(\alpha,\beta)}(x)  \nonumber \\
  & -(1-x)\bigl[\partial_x P_{m}^{(-\alpha,\beta)}(x)\bigr] P_n^{(\alpha,\beta)}(x)\Bigr\}, \qquad 
      n=0,1,2,\ldots,
\end{align*}
and
\begin{align*}
  &\hat{P}^{({\rm III},m)}_n(x;\alpha,\beta) \nonumber \\
  &\quad = \begin{cases}
     1 & \text{if $n=0$}, \\
     \frac{1}{\alpha+\beta+n-2m-1} \Bigl\{P_m^{(-\alpha,-\beta)}(x) [(x^2-1)\partial_x + (\alpha+\beta)x 
         +\alpha-\beta] P_{n-m-1}^{(\alpha,\beta)}(x) & \\
     \quad + (1-x^2) \bigl[\partial_x P_{m}^{(-\alpha,-\beta)}(x)\bigr] P_{n-m-1}^{(\alpha,\beta)}(x)\Bigr\}, 
         & \text{if $n=m+1$}, \\
      & \quad \text{$m+2,\ldots$}, \\  
  \end{cases}  
\end{align*}
after using (\ref{eq:A8}) and (\ref{eq:A9}).\par
%
%
It is worth observing that types I and II $X_m$-Jacobi EOPs are not independent since
\begin{equation*}
  \hat{P}^{({\rm II},m)}_n(x;\alpha,\beta) = (-1)^{n+m} \hat{P}^{({\rm I},m)}_n(-x;\beta,\alpha)
\end{equation*}
as a direct consequence of (\ref{eq:j1}), (\ref{eq:j2}), and the well-known symmetry relation $P^{(\alpha,\beta)}_n(x) = (-1)^n P^{(\beta,\alpha)}_n(-x)$ of Jacobi polynomials.\par
%
%
As we now plan to show, it is actually type II $X_m$-Jacobi EOPs (\ref{eq:j2}) that are considered in some mathematical studies \cite{bonneux, gomez13, liaw15}. Specifically, we are going to prove that $\hat{P}^{(\alpha,\beta)}_{m,n}(x)$, $n=m, m+1, m+2, \ldots$, defined in Eq.~(64) of \cite{gomez13} is such that
\begin{align}
  \hat{P}^{(\alpha,\beta)}_{m,n}(x) & = (-1)^m \frac{\alpha+1+n-2m}{\alpha+1+n-m} 
      \frac{\Gamma(2m-\alpha+\beta-1) \Gamma(2n-2m+\alpha+\beta+1)}{2^n m! (n-m)!
      \Gamma(m-\alpha+\beta-1) \Gamma(n-m+\alpha+\beta+1)} \nonumber \\
  & \quad \times  \hat{P}^{({\rm II},m)}_{n-m}(x;\alpha+1,\beta-1), \qquad n=m, m+1, m+2, \ldots.
      \label{eq:link}
\end{align} 
Equation (64) of \cite{gomez13}, rewritten in terms of monic Jacobi polynomials, is given by
\begin{align*}
  & \hat{P}^{(\alpha,\beta)}_{m,n}(x)  = \frac{(-1)^m}{\alpha+1+n-m} 
      \frac{\Gamma(2m-\alpha+\beta-1) \Gamma(2n-2m+\alpha+\beta+1)}{2^n m! (n-m)!
      \Gamma(m-\alpha+\beta-1) \Gamma(n-m+\alpha+\beta+1)} \\
  & \quad \times \Bigl[(n-m)(x-1) P^{(-\alpha-1,\beta-1)}_m(x) P^{(\alpha+2,\beta)}_{n-m-1}(x) + (\alpha+1-m)
      P^{(-\alpha-2,\beta)}_{m}(x) P^{(\alpha+1,\beta-1)}_{n-m}(x)\Bigr], \\
  & \qquad n=m, m+1, m+2, \ldots,
\end{align*}
while from (\ref{eq:j2}), we obtain
\begin{align*}
  & \hat{P}^{({\rm II},m)}_{n-m}(x;\alpha+1,\beta-1) = \frac{1}{\alpha+n+1-2m} \Bigl\{(n-m) (x-1)
       P^{(-\alpha-1,\beta-1)}_m(x)  P^{(\alpha+2,\beta)}_{n-m-1}(x) \\
  & \quad {}+ \Bigl[(\alpha+1) P^{(-\alpha-1,\beta-1)}_m(x) + m(1-x) P^{(-\alpha,\beta)}_{m-1}(x)\Bigr]
       P^{(\alpha+1,\beta-1)}_{n-m}(x)\Bigr\}, \\
  & \qquad n=m, m+1, m+2, \ldots.
\end{align*}
Hence, Eq.~(\ref{eq:link}) is equivalent to the following identity satisfied by Jacobi polynomials
\begin{equation}
  (\alpha+1) P^{(-\alpha-1,\beta-1)}_m(x) + m(1-x) P^{(-\alpha,\beta)}_{m-1}(x) = (\alpha+1-m)
  P^{(-\alpha-2,\beta}_m(x),  \label{eq:identity}
\end{equation}
which is proved in Appendix B.\par
%
%
\section{\boldmath LIMIT RELATIONS SATISFIED BY $X_m$ EOPS}

\setcounter{equation}{0}

The purpose of this Section is to extend to $X_m$ EOPs the well-known limit relations satisfied by Jacobi, Laguerre, and Hermite COPs \cite{koekoek}, which in monic form can be written as
\begin{equation}
  \lim_{\beta\to\infty} \beta^n P^{(\alpha,\beta)}_n\Bigl(1-\frac{2x}{\beta}\Bigr) = (-2)^n L^{(\alpha)}_n(x),
  \label{eq:jl}
\end{equation}
\begin{equation}
  \lim_{\alpha\to\infty} \alpha^{n/2} P^{(\alpha,\alpha)}_n\Bigl(\frac{x}{\sqrt{\alpha}}\Bigr) = H_n(x),
  \label{eq:jh}
\end{equation}
and
\begin{equation}
  \lim_{\alpha\to\infty} \frac{1}{(2\alpha)^{n/2}} L^{(\alpha)}_n \bigl(\sqrt{2\alpha}x+\alpha\bigr) =
   H_n(x).  \label{eq:lh}
\end{equation}
\par
%
%
\subsection{Going from Jacobi to Laguerre EOPs}

Equation (\ref{eq:jl}) can be generalized to $X_m$-Jacobi and $X_m$-Laguerre EOPs as follows:
\begin{equation}
  \lim_{\beta\to\infty} \beta^{n+m} \hat{P}^{({\rm I},m)}_n\left(1-\frac{2x}{\beta};\alpha,\beta\right) =
  (-2)^{n+m} \hat{L}^{({\rm I},m)}_n(x;\alpha), \qquad n=0,1,2,\ldots, \label{eq:jl1}
\end{equation}
\begin{equation}
  \lim_{\beta\to\infty} \beta^{n+m} \hat{P}^{({\rm II},m)}_n\left(1-\frac{2x}{\beta};\alpha,\beta\right) =
  (-2)^{n+m} \hat{L}^{({\rm II},m)}_n(x;\alpha), \qquad n=0,1,2,\ldots, \label{eq:jl2}
\end{equation}
and
\begin{equation}
  \lim_{\beta\to\infty} \beta^n \hat{P}^{({\rm III},m)}_n\left(1-\frac{2x}{\beta};\alpha,\beta\right) =
  (-2)^n \hat{L}^{({\rm III},m)}_n(x;\alpha), \qquad n=m+1,m+2,m+3,\ldots. \label{eq:jl3}
\end{equation}
\par
%
%
Let us first consider Eq.~(\ref{eq:jl1}). By using (\ref{eq:j1}), its left-hand side can be rewritten as
\begin{align*}
  & \lim_{\beta\to\infty} \beta^{n+m} \hat{P}^{({\rm I},m)}_n\Bigl(1-\frac{2x}{\beta};\alpha,\beta\Bigr) \\
  & \quad = \lim_{\beta\to\infty} \frac{\beta^{n+m}}{\beta+n-m} \Bigl\{P^{(\alpha,-\beta)}_m\Bigl(1-
       \frac{2x}{\beta}\Bigr) \Bigl[2n \Bigl(1-\frac{x}{\beta}\Bigr) P^{(\alpha+1,\beta+1)}_{n-1}\Bigl(1-
       \frac{2x}{\beta}\Bigr) \\
  & \qquad {}+ \beta P^{(\alpha,\beta)}_n\Bigl(1-\frac{2x}{\beta}\Bigr)\Bigr] - 2m\Bigl(1-\frac{x}{\beta}  
       \Bigr)
       P^{(\alpha+1,-\beta+1)}_{m-1}\Bigl(1-\frac{2x}{\beta}\Bigr) P^{(\alpha,\beta)}_n\Bigl(1-\frac{2x} 
       {\beta}\Bigr)\Bigr\} \\
  & \quad = \lim_{\beta\to\infty} \beta^{n+m-1} \Bigl\{P^{(\alpha,-\beta)}_m\Bigl(1-\frac{2x}{\beta}\Bigr)
       \Bigl[2n P^{(\alpha+1,\beta+1)}_{n-1}\Bigl(1-\frac{2x}{\beta}\Bigr) + \beta P^{(\alpha,\beta)}_n
       \Bigl(1-\frac{2x}{\beta}\Bigr)\Bigr] \\
  & \qquad {}- 2m P^{(\alpha+1,-\beta+1)}_{m-1}\Bigl(1-\frac{2x}{\beta}\Bigr) P^{(\alpha,\beta)}_n
       \Bigl(1-\frac{2x}{\beta}\Bigr)\Bigr\} \\
  & \quad = \Bigl[\lim_{\beta\to\infty} \beta^m P^{(\alpha,-\beta)}_m\Bigl(1-\frac{2x}{\beta}\Bigr)\Bigr]
       \Bigl[2n \lim_{\beta\to\infty} \beta^{n-1} P^{(\alpha+1,\beta+1)}_{n-1}\Bigl(1- \frac{2x}{\beta}\Bigr) 
       \\
  & \qquad {}+ \lim_{\beta\to\infty} \beta^n P^{(\alpha,\beta)}_n\Bigl(1-\frac{2x}{\beta}\Bigr)\Bigr]
       - 2m \Bigl[\lim_{\beta\to\infty} \beta^{m-1} P^{(\alpha+1,-\beta+1)}_{m-1}\Bigl(1-\frac{2x}{\beta}
       \Bigr)\Bigr] \\
  & \qquad \times \Bigl[\lim_{\beta\to\infty} \beta^n P^{(\alpha,\beta)}_n\Bigl(1-\frac{2x}{\beta}\Bigr)
       \Bigr].
\end{align*}
On employing now (\ref{eq:jl}), as well as its corollary
\begin{equation}
  \lim_{\beta\to\infty} \beta^n P^{(\alpha,-\beta)}_n\Bigl(1-\frac{2x}{\beta}\Bigr) = 2^n L^{(\alpha)}_n
  (-x),  \label{eq:jl-bis}
\end{equation}
we get
\begin{align*}
  & \lim_{\beta\to\infty} \beta^{n+m} \hat{P}^{({\rm I},m)}_n\Bigl(1-\frac{2x}{\beta};\alpha,\beta\Bigr) \\
  & \quad = 2^m L^{(\alpha)}_m(-x) \Bigl[2n (-2)^{n-1} L^{(\alpha+1)}_{n-1}(x) + (-2)^n L^{(\alpha)}_n
       (x)\Bigr] \\
  & \qquad {} - 2m 2^{m-1} L^{(\alpha+1)}_{m-1}(-x) (-2)^n L^{(\alpha)}_n(x) \\
  & \quad = (-2)^{n+m} (-1)^m \Bigl\{\Bigl[L^{(\alpha)}_m(-x) - m L^{(\alpha+1)}_{m-1}(-x)\Bigr] 
       L^{(\alpha)}_n(x) - n L^{(\alpha)}_m(-x) L^{(\alpha+1)}_{n-1}(x)\Bigr\} \\
  & \quad = (-2)^{n+m} (-1)^m \Bigl[L^{(\alpha+1)}_m(-x) L^{(\alpha)}_n(x) - n L^{(\alpha)}_m(-x)
       L^{(\alpha+1)}_{n-1}(x)\Bigr],
\end{align*}
where in the last step use is made of (\ref{eq:A4}). Comparison with (\ref{eq:l1}) completes the proof of (\ref{eq:jl1}).\par
%
%
Next, on inserting (\ref{eq:j2}) in the left-hand side of (\ref{eq:jl2}) and proceeding as in the previous case, we arrive at
\begin{align*}
  & \lim_{\beta\to\infty} \beta^{n+m} \hat{P}^{({\rm II},m)}_n\Bigl(1-\frac{2x}{\beta};\alpha,\beta\Bigr) \\
  & \quad = \frac{1}{m-n-\alpha} \Bigl\{2nx \Bigl[\lim_{\beta\to\infty} \beta^m P^{(-\alpha,\beta)}_m
      \Bigl(1-\frac{2x}{\beta}\bigr)\Bigl] \Bigl[\lim_{\beta\to\infty} \beta^{n-1}P^{(\alpha+1,\beta+1)}_{n-1}
      \Bigl(1- \frac{2x}{\beta}\Bigr)\Bigr] \\
  & \qquad {}-\alpha \Bigl[\lim_{\beta\to\infty} \beta^m P^{(-\alpha,\beta)}_m\Bigl(1-\frac{2x}{\beta}\bigr)
      \Bigr]\Bigl[\lim_{\beta\to\infty} \beta^n P^{(\alpha,\beta)}_n\Bigl(1-\frac{2x}{\beta}\Bigr)\Bigr] \\
  & \qquad {}-2mx \Bigl[\lim_{\beta\to\infty} \beta^{m-1} P^{(-\alpha+1,\beta+1)}_{m-1}\Bigl(1
      -\frac{2x}{\beta}\Bigr)\Bigr] \Bigl[\lim_{\beta\to\infty} \beta^n P^{(\alpha,\beta)}_n\Bigl(1
      -\frac{2x}{\beta}\Bigr)\Bigr]\Bigr\} \\
  & \quad = \frac{(-2)^{n+m}}{m-n-\alpha} \Bigl[-nx L^{(-\alpha)}_m(x)  L^{(\alpha+1)}_{n-1}(x)
      - \alpha L^{(-\alpha)}_m(x) L^{(\alpha)}_n(x) + mx L^{(-\alpha+1)}_{m-1}(x) L^{(\alpha)}_n(x)\Bigr] \\
  & \quad = \frac{(-2)^{n+m}}{\alpha+n-m} \Bigl\{nx L^{(-\alpha)}_m(x) L^{(\alpha+1)}_{n-1}(x)
      + \Bigl[\alpha L^{(-\alpha)}_m(x) - mx L^{(-\alpha+1)}_{m-1}(x)\Bigr] L^{(\alpha)}_n(x)\Bigr\}
\end{align*}
on using (\ref{eq:jl}). On comparing with (\ref{eq:l2}), to prove (\ref{eq:jl2}) it remains to check that
\begin{equation*}
  mx L^{(-\alpha+1)}_{m-1}(x) - \alpha L^{(-\alpha)}_m(x) = (m-\alpha) L^{(-\alpha-1)}_m(x).
\end{equation*}
From (\ref{eq:A4}), it results that such a relation is equivalent to
\begin{equation*}
  mx L^{(-\alpha+1)}_{m-1}(x) - \alpha L^{(-\alpha)}_m(x) = (m-\alpha) \Bigl[L^{(-\alpha)}_m(x) + m
  L^{(-\alpha)}_{m-1}(x)\Bigr]
\end{equation*}
or
\begin{equation*}
  x L^{(-\alpha+1)}_{m-1}(x) = (m-\alpha) L^{(-\alpha)}_{m-1}(x) + L^{(-\alpha)}_m(x),
\end{equation*}
which is indeed satisfied owing to (\ref{eq:A7}). This completes the proof of (\ref{eq:jl2}).\par
%
%
{}Finally, on inserting (\ref{eq:j3}) in the left-hand side of (\ref{eq:jl3}), we obtain that for $n=m+1,m+2,\ldots$
\begin{align*}
  & \lim_{\beta\to\infty} \beta^n \hat{P}^{({\rm III},m)}_n\Bigl(1-\frac{2x}{\beta};\alpha,\beta\Bigr) \\
  & \quad = \Bigl[\lim_{\beta\to\infty} \beta^m P^{(-\alpha	,-\beta)}_m\Bigl(1- \frac{2x}{\beta}\Bigr)
      \Bigr]\Bigl[\lim_{\beta\to\infty} \beta^{n-m} P^{(\alpha-1,\beta-1)}_{n-m}\Bigl(1-\frac{2x}{\beta}\Bigr)
      \Bigr] \\
  & \qquad {}+4mx \Bigl[\lim_{\beta\to\infty} \beta^{m-1} P^{(-\alpha+1,-\beta+1)}_{m-1}\Bigl(1
      -\frac{2x}{\beta}\Bigr)\Bigr] \Bigl[\lim_{\beta\to\infty} \beta^{n-m-1} P^{(\alpha,\beta)}_{n-m-1}
      \Bigl(1-\frac{2x}{\beta}\Bigr)\Bigr].
\end{align*}
On using (\ref{eq:jl}) and (\ref{eq:jl-bis}), we obtain
\begin{align*}
  & \lim_{\beta\to\infty} \beta^n \hat{P}^{({\rm III},m)}_n\Bigl(1-\frac{2x}{\beta};\alpha,\beta\Bigr) \\
  & \quad = (-1)^{n-m} 2^n \Bigl[L^{(-\alpha)}_m(-x) L^{(\alpha-1)}_{n-m}(x) - mx L^{(-\alpha+1)}_{m-1}
      (-x) L^{(\alpha)}_{n-m-1}(x)\Bigr].
\end{align*}
A comparison with (\ref{eq:l3}) then shows that Eq.~(\ref{eq:jl3}) is satisfied, which completes the proof.\par
%
%
\subsection{Going from Jacobi or Laguerre to Hermite EOPs }

For type III EOPs and even $m$ values, Eqs.~(\ref{eq:jh}) and (\ref{eq:lh}) can be generalized as follows:
\begin{equation}
  \lim_{\alpha\to\infty} \alpha^{n/2} \hat{P}^{({\rm III},m)}_n\Bigl(\frac{x}{\sqrt{\alpha}}; \alpha,\alpha
  \Bigr) = \hat{H}^{({\rm III},m)}_n(x), \qquad n=m+1,m+2,\ldots,  \label{eq:jh3}
\end{equation}
and
\begin{equation}
  \lim_{\alpha\to\infty} \frac{1}{(2\alpha)^{n/2}} \hat{L}^{({\rm III},m)}_n\bigl(\sqrt{2\alpha}x+\alpha;
  \alpha\bigr) = \hat{H}^{({\rm III},m)}_n(x), \qquad n=m+1,m+2,\ldots.  \label{eq:lh3}
\end{equation}
\par
%
%
Let us start with the proof of (\ref{eq:jh3}). From (\ref{eq:j3}), we obtain that for $n=m+1,m+2,\ldots$,
\begin{align*}
  &\lim_{\alpha\to\infty} \alpha^{n/2} \hat{P}^{({\rm III},m)}_n\Bigl(\frac{x}{\sqrt{\alpha}};\alpha,\alpha
       \Bigr) \\
  &\quad = \lim_{\alpha\to\infty} \frac{\alpha^{n/2}}{2\alpha+n-2m-1} \Bigl[(2\alpha+n-m-1) P^{(-\alpha,
      -\alpha)}_m\Bigl(\frac{x}{\sqrt{\alpha}}\Bigr) P^{(\alpha-1,\alpha-1)}_{n-m}\Bigl(\frac{x}{\sqrt{\alpha}}
      \Bigr) \\
  &\qquad{}+m\Bigl(1-\frac{x^2}{\alpha^2}\Bigr) P^{(-\alpha+1,-\alpha+1)}_{m-1}\Bigl(\frac{x}
      {\sqrt{\alpha}}\Bigr) P^{(\alpha,\alpha)}_{n-m-1}\Bigl(\frac{x}{\sqrt{\alpha}}\Bigr)\Bigr] \\
  & \quad = \Bigl[\lim_{\alpha\to\infty} \alpha^{m/2} P^{(-\alpha,-\alpha)}_m\Bigl(\frac{x}{\sqrt{\alpha}}
      \Bigr)\Bigr] \Bigl[\lim_{\alpha\to\infty} \alpha^{(n-m)/2} P^{(\alpha-1,\alpha-1)}_{n-m}\Bigl(
      \frac{x}{\sqrt{\alpha}}\Bigr)\Bigr] \\
  &\qquad{}+ \frac{m}{2} \Bigl[\lim_{\alpha\to\infty} \alpha^{(m-1)/2} P^{(-\alpha+1,-\alpha+1)}_{m-1}
      \Bigl(\frac{x}{\sqrt{\alpha}}\Bigr)\Bigr] \Bigl[\lim_{\alpha\to\infty} \alpha^{(n-m-1)/2} 
      P^{(\alpha,\alpha)}_{n-m-1}\Bigl(\frac{x}{\sqrt{\alpha}}\Bigr)\Bigr] \\
  &\quad = (-{\rm i})^m H_m({\rm i}x) H_{n-m}(x) + \frac{m}{2} (-{\rm i})^{m-1} H_{m-1}({\rm i}x)
      H_{n-m-1}(x) \\
  &\quad = \hat{H}^{({\rm III},m)}_n(x),
\end{align*}
where Eq.~(\ref{eq:jh}) and its corollary
\begin{equation*}
  \lim_{\alpha\to\infty} (-\alpha)^{n/2} P^{(-\alpha,-\alpha)}_n\Bigl(\frac{x}{\sqrt{\alpha}}\Bigr) =
  H_n({\rm i}x)
\end{equation*}
or
\begin{equation*}
  \lim_{\alpha\to\infty} \alpha^{n/2} P^{(-\alpha,-\alpha)}_n\Bigl(\frac{x}{\sqrt{\alpha}}\Bigr) =
  (-{\rm i})^n H_n({\rm i}x)
\end{equation*}
are used, as well as definition (\ref{eq:h}).\par
%
%
Let us now consider the proof of (\ref{eq:lh3}). This time, we start from Eq.~(\ref{eq:l3}) and obtain for $n=m+1,m+2,\ldots$,
\begin{align*}
  & \lim_{\alpha\to\infty} \frac{1}{(2\alpha)^{n/2}} \hat{L}^{({\rm III},m)}_n\bigl(\sqrt{2\alpha}x+\alpha;
       \alpha\bigr) \\
  &\quad = \lim_{\alpha\to\infty} \frac{1}{(2\alpha)^{n/2}} \Bigl[L^{(-\alpha)}_m\bigl(-\sqrt{2\alpha}x
      - \alpha\bigr) L^{(\alpha-1)}_{n-m}\bigl(\sqrt{2\alpha}x+\alpha\bigr) \\
  &\qquad{}- m\bigl(\sqrt{2\alpha}x+\alpha\bigr) L^{(-\alpha+1)}_{m-1}\bigl(-\sqrt{2\alpha}x-\alpha\bigr)
      L^{(\alpha)}_{n-m-1}\bigl(\sqrt{2\alpha}x+\alpha\bigr)\Bigr] \\
  &\quad = \Bigl[\lim_{\alpha\to\infty} \frac{1}{(2\alpha)^{m/2}} L^{(-\alpha)}_m\bigl(-\sqrt{2\alpha}
      x-\alpha\bigr)\Bigr] \Bigl[\lim_{\alpha\to\infty} \frac{1}{(2\alpha)^{(n-m)/2}} L^{(\alpha-1)}_{n-m}
      \bigl(\sqrt{2\alpha}x+\alpha\bigr)\Bigr] \\
  &\qquad{}- \frac{m}{2} \Bigl[\lim_{\alpha\to\infty} \frac{1}{(2\alpha)^{(m-1)/2}} L^{(-\alpha+1)}_{m-1}
      \bigl(-\sqrt{2\alpha}x-\alpha\bigr)\Bigr] \Bigl[\lim_{\alpha\to\infty} \frac{1}{(2\alpha)^{(n-m-1)/2}}
      L^{(\alpha)}_{n-m-1}\bigl(\sqrt{2\alpha}x+\alpha\bigr)\Bigr].
\end{align*}
Here, we employ (\ref{eq:lh}), as well as its corollary
\begin{equation*}
  \lim_{\alpha\to\infty} \frac{1}{(2\alpha)^{n/2}} L^{(-\alpha)}_n\bigl(-\sqrt{2\alpha}x-\alpha\bigr) =
  {\rm i}^n H_n({\rm i}x),
\end{equation*}
thus leading to
\begin{equation*}
  \lim_{\alpha\to\infty} \frac{1}{(2\alpha)^{n/2}} \hat{L}^{({\rm III},m)}_n\bigl(\sqrt{2\alpha}x+\alpha;
       \alpha\bigr) = {\rm i}^m H_m({\rm i}x) H_{n-m}(x) - \frac{m}{2} {\rm i}^{m-1} H_{m-1}({\rm i}x)
       H_{n-m-1}(x),
\end{equation*}
which amounts to $\hat{H}^{({\rm III},m)}_n(x)$, as given in (\ref{eq:h}).\par
%
%
In Appendix C, we provide some examples of the limit relations proved in this Section.\par
%
%
\section{QUADRATIC TRANSFORMATIONS RELATING HERMITE TO LAGUERRE EOPS}

\setcounter{equation}{0}

The purpose of this Section is to examine whether the quadratic transformations relating Hermite to Laguerre COPs \cite{koekoek}, which in monic form can be written as
\begin{equation}
  H_{2n}(x) = L^{(-1/2)}_n(x^2),  \qquad n=0, 1, 2, \ldots,   \label{eq:hl-even}
\end{equation}
and
\begin{equation}
  H_{2n+1}(x) = x L^{(1/2)}_n(x^2), \qquad n=0,1,2,\ldots,  \label{eq:hl-odd} 
\end{equation}
can be generalized to type III EOPs.\par
%
%
We wish to prove that the answer to this question is positive for even $X_{2m}$-Hermite EOPs and that in such a case the relation reads
\begin{equation}
  \hat{H}^{({\rm III},2m)}_{2n}(x) = \hat{L}^{({\rm III},m)}_n\Bigl(x^2;\frac{1}{2}\Bigr), \qquad n=0, m+1,
  m+2, \ldots.  \label{eq:hl3}
\end{equation}
\par
%
%
To start with, let us note that if we formally set $m=0$ and take Eqs.~(C.1) and (C.2) into account, Eq.~(\ref{eq:hl3}) reduces to the known relation (\ref{eq:hl-even}). It is also obvious that Eq.~(\ref{eq:hl3}) is fulfilled for $m=1, 2, 3, \ldots$ and $n=0$.\par
%
%
Let us therefore consider $m=1, 2, 3, \ldots$ and $n=m+1, m+2, \ldots$. Then from (\ref{eq:h}) it follows that
\begin{equation*}
  \hat{H}^{({\rm III}, 2m)}_{2n} (x) = (-1)^m [H_{2m}({\rm i}x)H_{2n-2m}(x) + {\rm i}m H_{2m-1}({\rm i}x)
  H_{2n-2m-1}(x)].
\end{equation*}
Next, Eqs.~(\ref{eq:hl-even}) and (\ref{eq:hl-odd}) imply that
\begin{equation*}
  \hat{H}^{({\rm III},2m)}_{2n}(x) = (-1)^m \Bigl[L^{(-1/2)}_m(-x^2) L^{(-1/2)}_{n-m}(x^2) - mx^2
  L^{(1/2)}_{m-1}(-x^2) L^{(1/2)}_{n-m-1}(x^2)\Bigr].  
\end{equation*}
 On comparing with
 \begin{equation}
  \hat{L}^{({\rm III},m)}_{n}(x^2;\alpha)= (-1)^m \Bigl[L^{(-\alpha)}_m(-x^2) L^{(\alpha-1)}_{n-m}(x^2)
  - mx^2 L^{(-\alpha+1)}_{m-1}(-x^2) L^{(\alpha)}_{n-m-1}(x^2)\Bigr] \label{eq:bracket}
\end{equation}
resulting from (\ref{eq:l3}), it follows that the right-hand sides of these relations coincide provided we set $\alpha=1/2$ in the second one. This completes the proof of (\ref{eq:hl3}).\par
%
%
Let us now turn ourselves to odd $X_{2m}$-Hermite EOPs, which from (\ref{eq:h}) can be written as
\begin{equation*}
  \hat{H}^{({\rm III},2m)}_{2n+1}(x) = (-1)^m [H_{2m}({\rm i}x) H_{2n-2m+1}(x) + {\rm i}m H_{2m-1}
  ({\rm i}x) H_{2n-2m}(x)]
\end{equation*}
for $n=m, m+1, m+2, \ldots$. On taking (\ref{eq:hl-even}) and (\ref{eq:hl-odd}) into account, this equation reduces to
\begin{equation*}
  \hat{H}^{({\rm III},2m)}_{2n+1}(x) = (-1)^m x \Bigl[L^{(-1/2)}_m(-x^2) L^{(1/2)}_{n-m}(x^2) - m
  L^{(1/2)}_{m-1}(-x^2) L^{(-1/2)}_{n-m}(x^2)\Bigr].
\end{equation*}
Here the square bracket on the right-hand side does not look to that in (\ref{eq:bracket}) for any value of $\alpha$.\par
%
%
In the special case where $n=m$, however, we directly get that
\begin{equation*}
  \hat{H}^{({\rm III},2m)}_{2m+1}(x) = (-1)^m [x H_{2m}({\rm i}x) + {\rm i}m H_{2m-1}({\rm i}x)],
\end{equation*}
which by a straightforward application of (\ref{eq:A1}) and (\ref{eq:A2}) reduces to
\begin{equation*}
  \hat{H}^{({\rm III},2m)}_{2m+1}(x) = (-1)^{m+1} {\rm i} H_{2m+1}({\rm i}x),
\end{equation*}
or, from (\ref{eq:hl-odd}),
\begin{equation*}
  \hat{H}^{({\rm III},2m)}_{2m+1}(x) = (-1)^m x L^{(1/2)}_m(-x^2).
\end{equation*}
This result may be compared with
\begin{equation*}
  \hat{L}^{({\rm III},m-1)}_m(x^2;\alpha) = (-1)^m L^{(-\alpha-1)}_m(-x^2),
\end{equation*}
which is a direct consequence of (\ref{eq:liaw}). Hence one may write
\begin{equation*}
  \hat{H}^{({\rm III},2m)}_{2m+1}(x)  = x \hat{L}^{({\rm III},m-1)}_m\Bigl(x^2;-\frac{3}{2}\Bigr), \qquad
   m=1, 2, 3, \ldots.
\end{equation*}
Apart from this result, the generalization of the Hermite-Laguerre quadratic transformation for odd $n$ values remains an open problem.\par
%
%
\section{CONCLUSION}

In the present paper, we have shown that the known asymptotic relations interconnecting Jacobi, Laguerre, and Hermite COPs can be generalized to the corresponding EOPs of codimension $m$.\par
%
%
For such a purpose, we have first listed the sets of EOPs to be considered, together with their explicit expressions in terms of COPs. We have also provided their link with some other definitions found in the literature.\par
%
%
We have then stated and proved limit relations allowing to get $X_m$-Laguerre EOPs of type I, II, or III from $X_m$-Jacobi EOPs of the same type, as well as $X_m$-Hermite EOPs of type III from $X_m$-Jacobi or $X_m$-Laguerre EOPs of type III.\par
%
%
{}Finally, we have established that a quadratic transformation applied to $X_m$-Laguerre EOPs of type III enables one to obtain even $X_{2m}$-Hermite EOPs. Whether a similar transformation would lead to odd $X_{2m}$-Hermite EOPs remains an open problem for future investigation.\par
%
%
Studying an extension of the limit relations presented here to more involved kinds of EOPs is also an open question for future work.\par
%
%
\section*{ACKNOWLEDGMENTS}

The author was supported by the Fonds de la Recherche Scientifique-FNRS under Grant No.~4.45.10.08.\par
%
%
\section*{AUTHOR DECLARATIONS}

\subsection*{Conflict of interest}

The author has no conflicts to disclose.\par
%
%
\subsection*{Author contributions}

C.\ Quesne: Conceptualization, Methodology, Investigation, Writing - review \& editing.\par
%
%
\section*{DATA AVAILABILITY}

Data sharing is not applicable to this article as no new data were created or analyzed in this study.\par
%
%
\section*{APPENDIX A: RELATIONS SATISFIED BY COPS}

\renewcommand{\theequation}{A.\arabic{equation}}
\setcounter{equation}{0}

In this Appendix, we list some relations satisfied by COPs \cite{gradshteyn, koekoek}, rewritten in monic form:
\begin{equation}
  \partial_x H_n(x) = n H_{n-1}(x),  \label{eq:A1}
\end{equation}
\begin{equation}
  (\partial_x - 2x) H_n(x) = -2 H_{n+1}(x),  \label{eq:A2}
\end{equation}
\begin{equation}
  \partial_x L^{(\alpha)}_n(x) = n L^{(\alpha+1)}_{n-1}(x),  \label{eq:A3}
\end{equation}
\begin{equation}
  L^{(\alpha-1)}_n(x) = L^{(\alpha)}_n(x) + n L^{(\alpha)}_{n-1}(x),  \label{eq:A4}
\end{equation}
\begin{equation}
  (x\partial_x + \alpha - x) L^{(\alpha)}_n(x) = - L^{(\alpha-1)}_{n+1}(x),  \label{eq:A5}
\end{equation}
\begin{equation}
  L^{(\alpha)}_{n+1}(x) + (2n+\alpha+1-x) L^{(\alpha)}_n(x) + n(n+\alpha) L^{(\alpha)}_{n-1}(x) = 0,
  \label{eq:A6}
\end{equation}
\begin{equation}
  x L^{(\alpha+1)}_n(x) = (n+\alpha+1) L^{(\alpha)}_n(x) + L^{(\alpha)}_{n+1}(x),  \label{eq:A7}
\end{equation}
\begin{equation}
  \partial_x P^{(\alpha,\beta)}_n(x) = n P^{(\alpha+1,\beta+1)}_{n-1}(x),  \label{eq:A8}
\end{equation}
\begin{equation}
  (1-x^2) \partial_x P^{(\alpha,\beta)}_n(x) + [\beta-\alpha-(\alpha+\beta)x] P^{(\alpha,\beta)}_n(x) =
  - (n+\alpha+\beta) P^{(\alpha-1,\beta-1)}_{n+1}(x),  \label{eq:A9}
\end{equation}
\begin{align}
  & (2n+\alpha+\beta+1) (2n+\alpha+\beta+2) \Bigl[(1-x) P^{(\alpha+1,\beta)}_n(x)
        + P^{(\alpha,\beta)}_{n+1}(x)\Bigr] \nonumber \\
  &\quad = 2(n+\alpha+1) (n+\alpha+\beta+1) P^{(\alpha,\beta)}_n(x),  \label{eq:A10}
\end{align}
\begin{equation}
  (2n+\alpha+\beta-1) (2n+\alpha+\beta) \Bigl[P^{(\alpha,\beta)}_n(x) - P^{(\alpha-1,\beta)}_n(x) \Bigr]
  = 2n(n+\beta) P^{(\alpha,\beta)}_{n-1}(x),  \label{eq:A11}
\end{equation}
\begin{equation}
  (2n+\alpha+\beta-1) \Bigl[P^{(\alpha,\beta-1)}_n(x) - P^{(\alpha-1,\beta)}_n(x)\Bigr] = 2n 
  P^{(\alpha,\beta)}_{n-1}(x).  \label{eq:A12}
\end{equation}
\par
%
%
\section*{APPENDIX B: PROOF OF AN IDENTITY SATISFIED BY JACOBI POLYNOMIALS}

\renewcommand{\theequation}{B.\arabic{equation}}
\setcounter{equation}{0}

The purpose of this Appendix is to prove Eq.~(\ref{eq:identity}).\par
%
%
On using Eq.~(\ref{eq:A10}) for $n\to m-1$, $\alpha\to -\alpha-1$, $\beta\to\beta$, we obtain
\begin{equation*}
  (1-x) P^{(-\alpha,\beta)}_{m-1}(x) = \frac{2(m-\alpha-1)(m-\alpha+\beta-1)}{(2m-\alpha+\beta-2)
  (2m-\alpha+\beta-1)} P^{(-\alpha-1,\beta)}_{m-1}(x) - P^{(-\alpha-1,\beta)}_m(x),
\end{equation*}
so the left-hand side of (\ref{eq:identity}) becomes
\begin{align*}
  & (\alpha+1) P^{(-\alpha-1,\beta-1)}_m(x) + m(1-x) P^{(-\alpha,\beta)}_{m-1}(x) = (\alpha+1) 
       P^{(-\alpha-1,\beta-1)}_m(x)\\
  &\quad {} + \frac{2m(m-\alpha-1)(m-\alpha+\beta-1)}
      {(2m-\alpha+\beta-2)(2m-\alpha+\beta-1)} P^{(-\alpha-1,\beta)}_{m-1}(x) - m P^{(-\alpha-1,\beta)}_m(x). 
\end{align*}
Equation (\ref{eq:A12}) for $n\to m$, $\alpha\to -\alpha-1$, $\beta\to\beta$, namely
\begin{equation*}
  P^{(-\alpha-1,\beta-1)}_m(x) = P^{(-\alpha-2,\beta)}_m(x) + \frac{2m}{2m-\alpha+\beta-2}
  P^{(-\alpha-1,\beta)}_{m-1}(x),
\end{equation*}
enables one to transform it into
\begin{align*}
  & (\alpha+1) P^{(-\alpha-1,\beta-1)}_m(x) + m(1-x) P^{(-\alpha,\beta)}_{m-1}(x) \\
  & \quad = (\alpha+1) P^{(-\alpha-2,\beta)}_m(x) + \frac{2m^2(m+\beta)}{(2m-\alpha+\beta-2)
      (2m-\alpha+\beta-1)} P^{(-\alpha-1,\beta)}_{m-1}(x) - m P^{(-\alpha-1,\beta)}_m(x).
\end{align*}
Finally, Eq.~(\ref{eq:A11}) for $n\to m$, $\alpha\to -\alpha-1$, $\beta\to\beta$ or
\begin{align*}
  & (2m-\alpha+\beta-2)(2m-\alpha+\beta-1) P^{(-\alpha-1,\beta)}_m(x) - 2m(m+\beta) P^{(-\alpha-1,\beta)}
  _{m-1}(x) \\
  & \quad = (2m-\alpha+\beta-2) (2m-\alpha+\beta-1) P^{(-\alpha-2,\beta)}_m(x)
\end{align*}
leads to the desired result (\ref{eq:identity}).\par
%
%
\section*{APPENDIX C: EXAMPLES OF LIMIT RELATIONS}

\renewcommand{\theequation}{C.\arabic{equation}}
\setcounter{equation}{0}

To start with, let us note that if we formally set $m=0$ in the EOPs definitions, we obtain that for $n=0, 1, 2, \ldots$,
\begin{equation}
  \hat{H}^{({\rm III},0)}_n(x) = H_n(x),  \label{eq:C1}
\end{equation}
\begin{equation}
  \hat{L}^{({\rm I},0)}_n(x;\alpha) = L^{(\alpha+1)}_n(x), \qquad \hat{L}^{({\rm II},0)}_n(x;\alpha) = 
  L^{(\alpha-1)}_n(x), \qquad \hat{L}^{({\rm III},0)}_n(x;\alpha) = L^{(\alpha-1)}_n(x), 
\end{equation}
\begin{align}
  & \hat{P}^{({\rm I},0)}_n(x;\alpha,\beta) = P^{(\alpha+1,\beta-1)}_n(x), \qquad 
       \hat{P}^{({\rm II},0)}_n(x;\alpha,\beta) = P^{(\alpha-1,\beta+1)}_n(x), \nonumber \\
  & \hat{P}^{({\rm III},0)}_n(x;\alpha,\beta) = P^{(\alpha-1,\beta+1)}_n(x),
\end{align}
after some straightforward application of identities satisfied by COPs. This shows that for $m=0$, the new limit relations of Section~3 connecting EOPs agree with the well-known ones satisfied by COPs.\par
%
%
For positive values of $m$, one gets for instance
\begin{equation*}
  \hat{H}^{({\rm III},2)}_3(x) = x^3 + \frac{3}{2}x,
\end{equation*}
\begin{equation*}
  \hat{L}^{({\rm I},1)}_2(x;\alpha) = x^3 - (\alpha+4)x^2 - (\alpha+1)(\alpha+4)x + (\alpha+1)(\alpha+2)
  (\alpha+4),
\end{equation*}
\begin{equation*}
  \hat{L}^{({\rm II},1)}_2(x;\alpha) = x^3 - (\alpha+2)x^2 - (\alpha-1)(\alpha+2)x + (\alpha-1)\alpha(\alpha+2),
\end{equation*}
\begin{equation*}
  \hat{L}^{({\rm III},2)}_3(x;\alpha) = x^3 - 3(\alpha-2)x^2 + 3(\alpha-1)(\alpha-2)x - \alpha(\alpha-1)(\alpha-2),
\end{equation*}
\begin{align*}
  &\hat{P}^{({\rm I},1)}_2(x;\alpha,\beta) = x^3 + \frac{2(\alpha-\beta+2)^2+(\alpha+\beta)(\alpha+\beta
      +6)}{(\alpha+\beta+4)(\alpha-\beta+2)}x^2 \\
  &\quad {} + \frac{(\alpha-\beta+2)^2+(\alpha+\beta)(2\alpha
      +2\beta+9)}{(\alpha+\beta+3)(\alpha+\beta+4)}x + \frac{(\alpha-\beta+2)^2(\alpha+\beta+2)-(\alpha
      +\beta)(\alpha+\beta+6)}{(\alpha-\beta+2)(\alpha+\beta+3)(\alpha+\beta+4)},
\end{align*}
\begin{align*}
  &\hat{P}^{({\rm II},1)}_2(x;\alpha, \beta) = x^3 + \frac{2(\alpha-\beta-2)^2+(\alpha+\beta)(\alpha+\beta
      +6)}{(\alpha+\beta+4)(\alpha-\beta-2)}x^2 \\
  & \quad {} + \frac{(\alpha-\beta-2)^2+(\alpha+\beta)(2\alpha
      +2\beta+9)}{(\alpha+\beta+3)(\alpha+\beta+4)}x +\frac{(\alpha-\beta-2)^2(\alpha+\beta+2)-(\alpha+
      \beta)(\alpha+\beta+6)}{(\alpha-\beta-2)(\alpha+\beta+3)(\alpha+\beta+4)},
\end{align*}
\begin{align*}
  &\hat{P}^{({\rm III},2)}_3(x;\alpha,\beta) = x^3 + \frac{3(\alpha-\beta)}{\alpha+\beta-4}x^2 + 3
      \frac{(\alpha-\beta)^2+\alpha+\beta-4}{(\alpha+\beta-3)(\alpha+\beta-4)}x \\
  & \quad {}+ \frac{(\alpha-\beta)[(\alpha-\beta)^2+3\alpha+3\beta-10]}{(\alpha+\beta-2)(\alpha
      +\beta-3)(\alpha+\beta-4)},
\end{align*}
on which one can directly check the relations
\begin{equation*}
  \lim_{\beta\to\infty} \beta^3 \hat{P}^{({\rm I},1)}_2 \Bigl(1 - \frac{2x}{\beta};\alpha,\beta\Bigr) =
  - 8 \hat{L}^{({\rm I},1)}_2(x;\alpha),
\end{equation*}
\begin{equation*}
  \lim_{\beta\to\infty} \beta^3 \hat{P}^{({\rm II},1)}_2 \Bigl(1 - \frac{2x}{\beta};\alpha,\beta\Bigr) =
  - 8 \hat{L}^{({\rm II},1)}_2(x;\alpha),
\end{equation*}
\begin{equation*}
  \lim_{\beta\to\infty} \beta^3 \hat{P}^{({\rm III},2)}_3 \Bigl(1 - \frac{2x}{\beta};\alpha,\beta\Bigr) =
  - 8 \hat{L}^{({\rm III},2)}_3(x;\alpha),
\end{equation*}
\begin{equation*}
  \lim_{\alpha\to\infty} \alpha^{3/2} \hat{P}^{({\rm III},2)}_3\Bigl(\frac{x}{\sqrt{\alpha}};\alpha,\alpha
  \Bigr) = \hat{H}^{({\rm III},2)}_3(x),
\end{equation*}
\begin{equation*}
  \lim_{\alpha\to\infty} \frac{1}{(2\alpha)^{3/2}} \hat{L}^{({\rm III},2)}_3\bigl(\sqrt{2\alpha}x+\alpha;
  \alpha\bigr) = \hat{H}^{({\rm III},2)}_3(x).
\end{equation*}
\par
%
%

\end{document}